\documentclass[12pt]{article}
\usepackage{amsthm, amsmath,amssymb,latexsym,color}

\newcommand{\beq}{\begin{equation}}
\newcommand{\eeq}{\end{equation}}

\newtheorem{dref}{Definition}[section] \newtheorem{lemma}[dref]{Lemma}
\newtheorem{theo}[dref]{Theorem} 
\newtheorem{remark}[dref]{Remark}

\newcommand{\ekv}[2]{\begin{equation}\label{#1}#2\end{equation}}

\def\supp{\mathop{\rm supp} \nolimits} 
\def\and {{\rm \; and \;}}

\def\Im {{\rm \, Im\,}}
\def\Re {{\rm \,Re\,}}

\newcommand {\pa}{\partial}
\newcommand {\ar}{\rightarrow}

\makeatletter
\@addtoreset{equation}{section}

\makeatother

\newtheorem{theorem}{Theorem} [section]

\newtheorem{proposition}[theorem]{Proposition}

\title{From resolvent bounds to semigroup bounds}
\author{B.~Helffer\footnote{Bernard.Helffer@math.u-psud.fr}\\Laboratoire de Mat\'ematiques, Univ Paris-Sud and
  CNRS,\\ F91405 Orsay Cedex France\\ \and\\ J.~Sj\"ostrand\footnote{Johannes.Sjostrand@u-bourgogne.fr}\\
IMB, UMR 5584, Universit\'e de Bourgogne,\\ 9 Av. A. Savary BP47870
F21078 Dijon Cedex France.}

\date{\today}
\begin{document}

\bibliographystyle{plain}

\maketitle
\begin{abstract}
The purpose of this note is to revisit the proof of the
Gearhardt-Pr\"uss-Hwang-Greiner theorem for a semigroup $S(t)$, following the
general idea of the proofs that we have seen in the literature and to
get an explicit estimate on $\Vert S(t)
\Vert$ in terms of bounds on the resolvent of the generator.
\end{abstract}
\section{Introduction}\label{int}

Let ${\cal H}$ be a complex Hilbert space and let $[0,+\infty [\ni
t\mapsto S(t)\in {\cal L}({\cal H},{\cal H})$ be a strongly continuous
semigroup with $S(0)=I$. Recall that by the Banach-Steinhaus theorem,
$\sup_J\Vert S(t)\Vert=:m(J)$ is bounded for every compact
interval $J\subset [0,+\infty [$. Using the semigroup property it
follows easily that there exist $M\ge 1$ and $\omega_0 \in \mathbb{R}$
such that $S(t)$ has the property 
\ekv{int.1}{P(M,\omega_0 ):\quad \Vert S(t)\Vert\le Me^{\omega_0 t},\ t\ge
0.}
In fact, we have this for $0\le t<1$ and for larger values of $t$,
write $t=[t]+r$, $[t]\in \mathbb{N} $, $0\le r<1$, and
$S(t)=S(1)^{[t]}S(r)$. 

Let $A$ be the generator of the semigroup (so that formally $S(t)=\exp tA$) and recall
(cf. \cite{EnNa07}, Chapter II or \cite{Paz}) that $A$ is closed and densely
defined. We also recall (\cite{EnNa07}, Theorem II.1.10) that 
\ekv{int.2}
{
(z-A)^{-1}=\int_0^\infty S(t)e^{-tz}dt,\quad \Vert (z-A)^{-1}\Vert \le
\frac{M}{\Re z-\omega_0 ,}
}
when $P(M,\omega_0 )$ holds and $z$ belongs to the open half-plane $\Re z> \omega_0 $. 

Recall the Hille-Yoshida theorem (\cite{EnNa07}, Th.~II.3.5) according
to which
the following three statements are equivalent when $\omega
\in \mathbb{R}$:
\begin{itemize}
\item $P(1,\omega )$ holds.
\item $\Vert (z-A)^{-1}\Vert \le (\Re z-\omega )^{-1}$, when $z\in
  \mathbb{C}$ and $\Re z> \omega $.
\item $\Vert (\lambda -A)^{-1}\Vert \le (\lambda -\omega )^{-1}$, when
  $\lambda \in ]\omega ,+\infty [$.
\end{itemize}
Here we may notice that we get from the special case $\omega =0$ to
general $\omega $ by passing from $S(t)$ to
$\widetilde{S}(t)=e^{-\omega t}S(t)$.

Also recall that there is a similar characterization of the property
$P(M,\omega )$ when $M>1$, in terms of the norms of all powers of the
resolvent. This is the Feller-Miyadera-Phillips theorem (\cite{EnNa07},
Th.~II.3.8). Since we need all powers of the resolvent, the practical
usefulness of that result is less evident. 

\par We next recall the Gearhardt-Pr\"uss-Hwang-Greiner theorem, see
\cite{EnNa07}, Theorem V.I.11, \cite{TrEm}, Theorem 19.1:
\begin{theo}\label{int1}~
\par\noindent 
(a) Assume that $\Vert (z-A)^{-1}\Vert$ is uniformly bounded in the
half-plane $\Re z\ge \omega $. Then there exists a constant $M>0$ such
that $P(M,\omega )$ holds.\\
(b) If $P(M,\omega )$ holds, then for every $\alpha >\omega $, $\Vert
(z-A)^{-1}\Vert$ is uniformly bounded in the half-plane $\Re z\ge
\alpha $.
\end{theo} 
The part (b) follows from (\ref{int.2}) with $\omega _0$ replaced by
$\omega $. 

{\bf The purpose of this note is to revisit the proof of (a), following the
general idea of the proofs that we have seen in the literature and to
get an explicit $t$ dependent estimate on $e^{-\omega t}\Vert S(t)
\Vert$, implying explicit bounds on $M$. }\\

 This idea is essentially to use that the resolvent and
  the inhomogeneous equation $(\partial _t-A)u=w$ in exponentially
  weighted spaces are related via
  Fourier-Laplace transform and we can use Plancherel's formula. Variants of
  this simple idea have also been used in more concrete
  situations. See \cite{BuZw, GGN, Hi, Sch09}.\\

Note that we can improve a little the conclusion of (a). If the
property (a) is true for some $\omega$ then it is automatically true
 for some  $\omega '<\omega $. We recall indeed 
  the following
\begin{lemma}\label{int01}~\\ If for some $r(\omega )>0$, $\Vert (z-A)^{-1}\Vert \le \frac{1}{r(\omega )}$
  for $\Re z>\omega $, then for every $\omega '\in ]\omega
  -r(\omega ),\omega ]$ we have 
$$
\Vert (z-A)^{-1}\Vert \le \frac{1}{r(\omega )-(\omega -\omega
  ')},\ \Re z>\omega '.
$$
\end{lemma} 
\medskip
\par\noindent {\bf Proof.} Let $\widetilde{z}\in \mathbb C $, $\Re
\widetilde{z}>\omega $. Then $\Vert (\widetilde{z}-A)^{-1}\Vert \le
\frac{1}{r(\omega )}$. For $z\in \mathbb C $ with $|z-\widetilde{z}| <r(\omega
)$, we have 
$$
(z-A)(\widetilde{z}-A)^{-1}=1+(z-\widetilde{z})(\widetilde{z}-A)^{-1},\hbox{
where }\Vert (z-\widetilde{z})(\widetilde{z}-A)^{-1}\Vert \le
|z-\widetilde{z}|/r(\omega )<1,
$$ 
so $1+(z-\widetilde{z})(\widetilde{z}-A)^{-1}$ is invertible and 
$$
\Vert (1+(z-\widetilde{z})(\widetilde{z}-A)^{-1})^{-1}\Vert \le
\frac{1}{1-|z-\widetilde{z}|/r(\omega )}.
$$
Hence $z$ belongs to the resolvent set of $A$ and 
$$
(z-A)^{-1}=(\widetilde{z}-A)^{-1}(1+(z-\widetilde{z})(\widetilde{z}-A)^{-1})^{-1},\
\Vert (z-A)^{-1}\Vert \le \frac{1}{r(\omega )-|z-\widetilde{z}|}.
$$

Now, if $z\in \mathbb C$ and $\Re z>\omega '$, we can find
$\widetilde{z}\in \mathbb C$ with $\Re z>\omega $,
$|z-\widetilde{z}|<\omega -\omega '$ and the lemma follows.\hfill{$\Box$}

\begin{remark}\label{Rem1.3}~\\
Let
$$\omega _0=\inf\{\omega \in \mathbb R \, \{z\in \mathbb C ; \Re z>\omega
\}\subset \rho (A) \hbox{ and }\sup_{\Re z>\omega
}\Vert (z-A)^{-1}\Vert <\infty \}.$$
For $\omega >\omega _0$, we may define $r(\omega )$ by
$$
\frac{1}{r(\omega )}=\sup_{\Re z>\omega }\Vert (z-A)^{-1}\Vert.
$$
Then $r(\omega )$ is an increasing function of $\omega $;  for every $\omega \in ]\omega _0,\infty [$, we have $\omega
-r(\omega )\ge \omega _0$ and for $\omega '\in [\omega -r(\omega
),\omega ]$ we have 
$$
r(\omega ')\ge r(\omega )-(\omega -\omega ').
$$
We may state all this more elegantly by saying that $r$ is a Lipschitz
function on $]\omega _0,+\infty [$ satisfying
$$
0\le \frac{d r}{d\omega }\le 1\,.
$$
 Moreover, if $\omega _0>-\infty $, then $r(\omega )\to 0$ when
$\omega \searrow \omega _0$.\end{remark}

\begin{remark}~\\
 Notice that by (\ref{int.1}), (\ref{int.2}), we already know that
$\Vert (z-A)^{-1}\Vert$ is uniformly bounded in the half-plane $\Re
z\ge \beta $, if $\beta >\omega _0$. If $\alpha \le \omega _0$, we
see that $\Vert (z-A)^{-1}\Vert$ is uniformly bounded in the
half-plane $\Re z\ge \alpha $, provided that 
\begin{itemize}
\item we have this uniform boundedness on the line $\Re z=\alpha $,
\item $A$ has no spectrum in the half-plane $\Re z\ge \alpha $,
\item $\Vert (z-A)^{-1}\Vert$ does not grow too wildly in the strip
  $\alpha \le \Re z \le \beta $: $\Vert (z-A)^{-1}\Vert \le {\cal
    O}(1) \exp ({\cal O}(1)\exp (k|\Im z|))$, where $k<\pi /(\beta
  -\alpha )$.

\end{itemize}
We then also have
\ekv{int.3}{
\sup_{\Re z\ge \alpha }\Vert (z-A)^{-1}\Vert = \sup_{\Re z=\alpha
}\Vert (z-A)^{-1}\Vert .} 
This follows from the subharmonicity of $\ln ||(z-A)^{-1}||$,
Hadamard's theorem (or Phragm\'en-Lindel\"of in
exponential coordinates) and the maximum principle. 
\end{remark}

\par Our main result is:
\begin{theo}\label{int2}~\\
We make the assumptions of Theorem \ref{int1}, (a) and define
$r(\omega )>0$ by
$$
\frac{1}{r(\omega )}=\sup_{\Re z\ge \omega }\Vert (z-A)^{-1}\Vert.
$$
Let $m(t)\ge \Vert S(t)\Vert$ be a continuous positive function.
 Then
for all $t,a,\widetilde{a}>0$, such that $t=a+\widetilde{a}$, we have
\ekv{int.4} { \Vert S(t)\Vert \le \frac{e^{\omega t}} {r(\omega )\Vert
    \frac{1}{m}\Vert_{e^{-\omega \cdot }L^2([0,a])}\Vert
    \frac{1}{m}\Vert_{e^{-\omega \cdot }L^2([0,\widetilde{a}])}}.}
\end{theo}
Here the norms are always the natural ones obtained from ${\cal H}$,
$L^2$, thus for instance 
$\Vert S(t)\Vert= \Vert S(t)\Vert_{{\cal L}({\cal H},{\cal H})}$, if
$u$ is a function on $\mathbb{R}$ with values in $\mathbb{C}$ or in
${\cal H}$, $\Vert u\Vert$  denotes the natural $L^2$ norm, when
the norm is taken over a subset $J$ of $\mathbb{R}$, this is indicated
with a ``$L^2(J)$''. In (\ref{int.4}) we also have the natural norm
in the exponentially weighted space $e^{-\omega \cdot }L^2([0,a])$ and
similarly with $\widetilde{a}$ instead of $a$; $\Vert
f\Vert_{e^{-\omega \cdot }L^2([0,a])}=\Vert e^{\omega \cdot }f(\cdot
)\Vert_{L^2([0,a])}$. \\

As we shall see in the next section, under the assumption of the
theorem, we have $P(M,\omega)$
 with an explicit $M$. See also the appendix.

We also have the following variant of the main result that can be
useful in problems of return to equilibrium.

\begin{theo}\label{int3}~\\
  We make the assumptions of Theorem \ref{int2}, so that (\ref{int.4})
  holds. Let $\widetilde{\omega }<\omega $ and assume that $A$ has no
  spectrum on the line $\Re z=\widetilde{\omega }$ and that the
  spectrum of $A$ in the half-plane $\Re z>\widetilde{\omega }$ is
  compact (and included in the strip $\widetilde{\omega }<\Re z<\omega
  $). Assume that $\Vert (z-A)^{-1}\Vert$ is uniformly bounded on $\{
  z\in \mathbb{C};\, \Re z\ge \widetilde{\omega }\} \setminus U$,
  where $U$ is any neighborhood of $\sigma _+(A):=\{ z\in \sigma
  (A);\, \Re z>\widetilde{\omega }\} $ and define $r(\widetilde{\omega
  })$ by
$$
\frac{1}{r(\widetilde{\omega })}=\sup_{\Re z=\widetilde{\omega }}\Vert
(z-A)^{-1}\Vert .
$$ 
Then for every $t>0$,
$$
S(t)=S(t)\Pi _++R(t)=S(t)\Pi _++S(t)(1-\Pi _+),
$$
where for all $a,\widetilde{a}>0$ with $a+\widetilde{a}=t$,
\ekv{int.5} { \Vert R(t)\Vert \le
  \frac{e^{\widetilde{\omega} t}}
{r(\widetilde{\omega} )\Vert \frac{1}{m}\Vert_{e^{-\widetilde{\omega} \cdot
      }L^2([0,a])}\Vert \frac{1}{m}\Vert_{e^{-\widetilde{\omega} \cdot
      }L^2([0,\widetilde{a}])}}\Vert I-\Pi _+\Vert.}
Here $\Pi _+$ denotes the spectral projection associated to $\sigma
_+(A)$: 
$$
\Pi _+=\frac{1}{2\pi i}\int_{\partial V}(z-A)^{-1}dz,
$$
where $V$ is any compact neighborhood of $\sigma _+(A)$ with $C^1$
boundary, disjoint from $\sigma (A)\setminus \sigma _+(A)$.   
\end{theo}

\section{Applications : Explicit bounds in the abstract framework}
Theorem \ref{int2} has two ingredients: the existence of some initial
control
 by $m(t)$ and the additional information on the resolvent.
\subsection{A quantitative Gearhardt-Pr\"uss statement}\label{ss2.1}
 As
 observed 
 in the introduction (see \eqref{int.1}), we have at least an estimate  with 
 $m(t) = \widehat M \,  \exp \widehat \omega t$,
 for some $\widehat \omega \geq \omega$. We apply Theorem \ref{int2}
 with this $m(t)$ and   $a=\tilde a = \frac{t}{2}$. The term appearing
 in the denominator of \eqref{int.4} becomes  
\begin{equation}
\Vert \frac{1}{m}\Vert_{e^{-\omega \cdot
      }L^2([0,a])}\Vert \frac{1}{m}\Vert_{e^{-\omega \cdot
      }L^2([0,\widetilde{a}])}
 =\frac 12 \, {\widehat M}^{-2} \, t \,,
\end{equation}
if $\widehat \omega=\omega$, and
\begin{equation}
=  \frac{1}{2 \widehat{M}^{2} (\widehat \omega-\omega)} \left[1 - \exp ((
  \omega-\widehat \omega) t)\right]\,,
\end{equation}
if $\widehat \omega > \omega $.

\par Hence we obtain the estimate with a new  $m^{new}(t)$, with
$$
m^{new}(t)= \frac{ 2 \widehat{M}^{2} (\widehat
  \omega-\omega)}{r(\omega)  [1 - \exp ((
  \omega-\widehat \omega) t)] } \;    \exp \omega t .
$$
This gives in particular that $S(t)$ satisfies $P( M,\omega)$,
 with
$$
M = \sup_t \left(\exp - \omega t\, \min (\widehat M \exp \widehat
\omega t ,  m^{new}(t))\right)\,.
$$
We will see how to optimize over $\omega$ in Subsection \ref{ss2.3}.
 
Let us push the computation. Without loss of generality,
we can assume $\widehat \omega=0$ and we make the assumption
in Theorem \ref{int2} for some $\omega <0$.  Combining 
 Theorem \ref{int2}  and the trivial estimate
$$ 
|| S(t)|| \leq  \widehat M = \widehat M  \exp - \omega t\, \exp \omega t
$$
we obtain that we have $P(M,\omega)$ with 
$$
M = \widehat M  \sup_t \left( \min (\exp - \omega t , \frac{2 \widehat
  M  |\omega|}{r(\omega) (1- \exp  \omega t
)})\right).
$$
This can be rewritten in the form:
$$
M = \widehat M  \sup_{u\in ]0,1[}\left( \min ( \frac{1}{u},  \frac{ 2
  \widehat M |\omega|}
{r(\omega) (1-u)}\right)  = 1 + 2 \frac{\widehat M |\omega|}{ r(\omega)}\,.
$$
\begin{proposition}\label{propa}~\\
Let $S(t)$ be a  continuous semigroup  such that $P(\widehat
M,\widehat \omega)$
 is satisfied for some pair $(\widehat M, \widehat \omega)$  and such that  $r(\omega) >0$
for some $\omega  <  \widehat \omega $. Then:
\begin{equation}\label{contra}
|| S(t)|| \leq \widehat M  \left(1 +  \frac{2 \widehat M  (
  \widehat \omega -\omega) }{ r(\omega)}\right)\; \exp \omega t\,.
\end{equation}
\end{proposition}

\subsection{Estimate with exponential gain.}\label{ss2.2}

In the same spirit, and combining with Lemma \ref{int01}, we get
 the following extension of 
\eqref{contra}  (with $\widehat \omega=0$) 
  \begin{equation}\label{contrb}
|| S(t)|| \leq  \widehat M \, \left(\frac{(1-s) r(\omega) + 2 \widehat
M ( \widehat \omega -\omega  +  s r(\omega))}{(1-s)r(\omega )}\right)    \;   \exp (\omega - s r(\omega)) t\,,\,
\forall s\in [0,1[\,.
\end{equation}
 Taking $s= \frac{t}{1 +t}$  gives a rather
  optimal decay
 at $\infty$ in $\mathcal O (t) \exp (\omega
  -r(\omega) ) t$.\\

If we assume now  instead  the control of the norm of the resolvent
  on  $\Re z \geq 0$, hence if we are in the case $\omega =\widehat
  \omega =0$, we get 
$$
|| S (t)|| \leq \frac{2 \widehat M }{r(0) t}\,,
$$
 and using the semi-group property $\leq \left(\frac{2 \widehat M N}{r(0)
   t}\right)^N$, for any $N\geq 1$.
Hence we can get  an explicit control of the decay of $S(t)$, by
optimizing over $N$. As in the theory of analytic symbols, we can take
 $N= E(\alpha t)$ where $E(s)$ denotes the integer part of $s$ and
$\alpha$
 such that $\alpha < r(0) / (2 \widehat M)$, we get an exponential
 decay of $S(t)$.\\
Alternately, we can use the extension of the resolvent on $\Re z > -s
r(0)$
 and this leads to~:
 \begin{equation}\label{contrbprime}
|| S(t)|| \leq  \widehat M \, \left(\frac{(1-s) + 2 \widehat
M s }{(1-s)}\right)    \;   \exp (- s r(0)) t\,,\,
\forall s\in [0,1[\,.
\end{equation}

\subsection{The limit $\omega \searrow \omega_0$}\label{ss2.3}~\\
Consider the situation of Theorem \ref{int2} and let $\omega _0$ be as
in Remark \ref{Rem1.3}. Assume that $\omega _0>-\infty $  so that $r(\omega )\to 0$, when
$\omega \to \omega _0$.  For $t\ge 1$, $\omega >\omega _0$, we
get from (\ref{int.4}):
\begin{equation}\label{rem.1new}
e^{-\omega _0t}\Vert S(t)\Vert \le \frac{e^{t(\omega -\omega
    _0)}}{r(\omega )\int_0^{1/2}m(s)^{-2}e^{2\omega _0s}ds}\le {\cal
  O}(1)\frac{e^{t(\omega -\omega _0)}}{r(\omega )}.
\end{equation} 
Optimizing over $\omega \in ]\omega_0,\omega_0+\epsilon_0]$, we  get
 the existence of $C$ such that
\begin{equation}
e^{-\omega _0t}\Vert S(t)\Vert \le C \exp \Phi(t)\,,
\end{equation}
with
$$
\Phi(t)=\inf_{\omega \in ]\omega_0,\omega_0+\epsilon_0]} t(\omega -
  \omega_0) - \ln r(\omega)\,.
$$
It is clear that $\lim_{t\ar +\infty} \Phi(t)/t=0$, but to have
 a more quantitative version, we need  some information
 on the behavior of $r(\omega)$ as $\omega\searrow \omega_0$. Let us
 treat two examples.\\

 If 
$$
r(\omega )\ge \frac{(\omega -\omega _0)^k}{C},\hbox{ when }0<\omega
-\omega _0\ll 1,
$$
for some constants $C,k>0$, then choosing $\omega -\omega _0=k/t$ in
(\ref{rem.1new}), we get 
$$
e^{-\omega _0t}\Vert S(t)\Vert \le {\cal O}(1)t^k,\ t\ge 1.
$$

\par On the other hand, if 
$$
r(\omega )\ge \exp -\frac{(\omega -\omega _0)^{-\alpha }}{C\alpha
},\hbox{ when }0<\omega -\omega _0\ll 1,
$$
for some constants $C,\alpha >0$, then 
$$
\frac{e^{t(\omega -\omega _0)}}{r(\omega )}\le \exp \left( t(\omega
  -\omega _0)+\frac{(\omega -\omega _0)^{-\alpha }}{C\alpha }  \right),
$$
and choosing $\omega -\omega _0=(Ct)^{-\frac{1}{\alpha +1}}$, we get 
 the existence of a constant $\widehat C$ such that
$$
e^{-\omega _0t}\Vert S(t)\Vert\leq e^{\widehat C t^{\frac{\alpha
  }{\alpha +1}}},\ t\ge 1.
$$

\section{Applications to concrete examples}

\subsection{The  complex Airy operator on the half-line}

Let us consider (as in \cite{Alm}) the Dirichlet realization $P^D$  of the 
Airy operator on $\mathbb R^+$~:  $D_x^2 + ix$ and $P$ the realization
of $D_x^2+ ix$ in $\mathbb R$. 
One can determine explicitly its spectrum (using Sibuya's theory or
Combes-Thomas's trick) as 
$$
\sigma (P^D):=\{  \lambda_j \,  e^{ i \frac \pi 3}  \,,\, j\in \mathbb
N^*\}
$$
where the $\lambda_j$'s are the eigenvalues
 (immediately related to the zeroes of the Airy function)
 of the Dirichlet realization in $\mathbb R^+$ of $D_x^2 +x$.\\

It was shown in  \cite{He}, that
 $|| (P^D -z)^{-1}||$
is  as $\Re z >0$ and $\Im z \mapsto +\infty$ 
 asymptotically equivalent to $|| (P -\Re z)^{-1}||$
 and  that  $|| (P^D -z)^{-1}||$ tends to $0$
  as $\Re z >0$ and $\Im z \mapsto -\infty$. The standard
 Gearhardt-Pr\"uss theorem, applied to $A:=-P^D$, permits to show that,
 for any $\omega >  -\lambda_1\,  \cos \frac \pi 3 \,$, we have
$$
|| S(t)||\leq M_\omega \exp (\omega t)\,.
$$
Theorem \ref{int3} permits the following improvment~:
$$
S (t) = \exp \left(-e^{ i \frac{\pi}{3}}\, \lambda_1\,t \right) \; \Pi_+ + R(t)\,,
$$
with
$$
|| R(t) ||\leq M_{\tilde \omega} \exp (\tilde \omega t)\,,
$$
for any  $\tilde \omega >  -  \lambda_2\, \cos \frac \pi 3 \,$.\\
Here $\Pi_+$ is the projector associated with the eigenfunction of
$P^D$
 associated with $ \lambda_1\, e^{ i\frac \pi 3}\,$.
Hence we get a much better control of the semi-group.

 \subsection{The case of the Kramers-Fokker-Planck operator}

Inspired by the work by F.~H\'erau and F.~Nier
\cite{HeNi04}, F.~H\'erau, J.~Sj\"ostrand and C.~Stolk \cite{HeSjSt05}
studied the Kramers-Fokker-Planck operator \ekv{int.6} { P=y\cdot
  h\partial _x-V'(x)\cdot h\partial _x+\frac{\gamma }{2}(y-h\partial
  _y)(y+h\partial _y) } on $\mathbb{R}^{2n}=\mathbb{R}_x^n\times
\mathbb{R}_y^n $, where $\gamma >0$ is fixed and we let $h\to 0$. We
assume that $V\in C^\infty (\mathbb{R}^n;\mathbb{R})$ with $\partial
^\alpha V={\cal O}(1)$ for every $\alpha \in \mathbb{N}^n$ of length
$\ge 2$ and we also assume that $V$ is a Morse function such that
$|\nabla V(x)|\ge 1/C$ when $|x|\ge C$ for some constant $C>0$. Then
we know from \cite{HeNi04} and under much weaker assumptions from
B.~Helffer, F.~Nier \cite{HelNi} that $P$ is maximally accretive with
$\Re P\ge 0$, so that $P$ generates a semi-group of contractions
$e^{-tP/h}$, $t\ge 0$. In particular the spectrum of $P$ is contained
in the closed right half plane. In \cite{HeSjSt05} it was shown that
for every fixed $C>0$ and for $h>0$ small enough, the spectrum of $P$
in the strip $0\le \Re z \le Ch$ is discrete and the eigenvalues are
of the form \ekv{int.7} { E_j=\lambda _jh+o(h),\quad \Re \lambda_j\le Ch, } where $\lambda _j$
are eigenvalues of the different quadratic approximations of $P_{h=1}$
at the various points $(x_k,0)$ where $V'(x_k)=0$. Here the points
$E_j$ all belong to a sector $|\Im \lambda|\le {\cal O}(\Re \lambda)$, so the
eigenvalues in (\ref{int.7}) are all confined to a disc
$D(0,\widetilde{C}h)$.

\par It was also shown in \cite{HeSjSt05} that if $\widetilde{\omega
}\ge 0$ and $\Re \lambda _j\ne \widetilde{\omega }$ for all the
eigenvalues $\lambda _j$, then $\Vert (P-z)^{-1}\Vert ={\cal O}(1/h)$
uniformly on the line $\Re z=h\widetilde{\omega }$. The same estimate
holds when $0\le \Re z\le Ch$ and $|z|\ge \widetilde{C}h$. Actually, using a
form of semi-classical sub-ellipticity (closely related in spirit to
the one established in \cite{HeNi04} and further studied in \cite{HelNi}) it was also shown that this
estimate holds in a larger parabolic neighborhood of $i\mathbb{R}$
away from the disc $D(0,\widetilde{C}h)$, and using this stronger
result and a contour deformation in a standard integral representation
of $e^{-tP/h}$ (again in the spirit of
\cite{HeNi04}) it was established in \cite{HeSjSt05} that \ekv{int.8}
{ e^{-tP/h}=e^{-tP/h}\Pi _++R(t), } where $\Pi _+$ is the spectral
projection associated with  $\{ z\in \sigma (P);\, 0\le \Re z\le
\widetilde{\omega }\}$, and $\Vert R(t)\Vert \le
\mathrm{Const.\,}e^{-t\widetilde{\omega }}$. Now this result becomes
a direct application of Theorem \ref{int3} to $A:=-P/h$ and we do not
need any bounds on the resolvent in the region $\Re
z>h\widetilde{\omega }$.

\par In \cite{HeHiSj08a, HeHiSj08b} similar results were obtained for
more general operators, for which we do not necessarily have any bound on
the resolvent beyond a strip, and the proof was to use microlocal
coercivity outside a compact set in slightly weighted
$L^2$-spaces. Again Theorem \ref{int3} would give some
simplifications. 

\subsection{The complex harmonic oscillator}
The complex harmonic oscillator
$$
P:= D_x^2 + i x^2$$ on the line  was studied by E.B.~Davies
\cite{Dav2, Dav3}, L.~Boulton, \cite{Bou02}
 and M.~Zworski \cite{Zw} in connection with the analysis of the
 pseudospectra.  As for the complex Airy operator, it is easy to determine
 the spectrum which is given by $e^{i \frac \pi 4} (2j+1)\,,\, j\in
 \mathbb N$. This operator is maximally accretive
 and we can apply Theorem \ref{int3} with $A = - P$.
From these works as well as those of K.~Pravda Starov \cite{Pr06} and
Dencker-Sj\"ostrand-Zworski \cite{DSZ}, we know that for fixed $\Re z$ as
 $\Im z \ar +\infty$, 
$$ 
\lim_{\Im z \ar + \infty} || (P-z)^{-1}|| =0\,.
$$

\par More precisely, for any compact interval $K$, there exists $C>0$ such
that
$$
|| (P-z)^{-1}||  \leq C \, |\Im z|^{-\frac 1 3}\,,\, \mbox{ for } \Im z
\geq C, \Re z \in K \,.
$$
This follows from \cite{Pr06, DSZ}, notice here that the
results in \cite{DSZ} are given in the semi-classical limit for the
spectral parameter in a compact set, but there is a simple scaling
argument, allowing to pass to the limit of high frequency. See for
example \cite{Sj, Sj2}. As $\Im z \ar -\infty$ we have by more elementary estimates:
$$
|| (P-z)^{-1}||  \leq \, |\Im z|^{-1}\,,\, \mbox{ for } \Im z
<0\,.
$$

We can therefore apply Theorem \ref{int3} and get
$$
S (t) = \exp \left(-e^{ i \frac{\pi}{4}} \,t \right) \; \Pi_+ + R(t)\,,
$$
with
$$
|| R(t) ||\leq M_{\tilde \omega} \exp (\tilde \omega t)\,,
$$
for any  $\tilde \omega >  - 3 \, \cos \frac \pi 4 \,$.
Here $\Pi_+$ is the spectral projection associated with the
eigenvalue $e^{ i\frac \pi 4}$ of
$P$.

Hence we get again  a much better control of the semi-group.

\section{Proofs of the main statements}
\subsection{Proof of Theorem \ref{int2}}\label{pr2}
\setcounter{equation}{0}

As already mentioned, we shall use the inhomogeneous equation
\begin{equation}\label{3}
(\partial _t -A) u = w  \hbox{ on }\mathbb{R}.
\end{equation}
Recall that if $v\in {\cal H}$, then $S(t)v\in C^0([0,\infty [;{\cal
  H})$, while if $v\in {\cal D}(A)$, then $S(t)v\in C^1([0,\infty [;{\cal
  H})\cap C^0([0,\infty [;{\cal
  D}(A))$ and 
\ekv{ny.1}
{
AS(t)v=S(t)Av,\quad (\partial _t-A)S(t)v=0\,.
}
\par Let $C_+^0 ({\cal H})$ denote the subspace of all $v\in C^0
(\mathbb{R};{\cal H})$ that vanish near $-\infty $. For $k\in
\mathbb{N}$, we define $C_+^k({\cal H})$ and $C_+^k({\cal D}(A))$
similarly. For $w\in C_+^0({\cal H})$, we define $Ew \in C_+^0({\cal
  H})$ by \ekv{ny.2} { Ew(t)=\int_{-\infty }^t S(t-s)w(s)ds.  } It is
easy to see that $E$ is continuous: $C_+^k({\cal H})\to C_+^k({\cal
  H})$, $C_+^k({\cal D}(A))\to C_+^k({\cal D}(A))$ and if $w\in
C_+^1({\cal H})\cap C_+^0({\cal D}(A))$, then $u=Ew$ is the unique
solution in the same space of (\ref{3}). More precisely, we have
\ekv{ny.3} { (\partial _t-A)Ew=w,\quad E(\partial _t-A)u=u, } for all
$u,w\in C_+^1({\cal H})\cap C_+^0({\cal D}(A))$

Now recall that we have $P(M,\omega _0)$ in (\ref{int.1}) for some
$M,\omega _0$. If $\omega _1>\omega _0$ and $w\in C_+^0({\cal H})\cap
e^{\omega _1\cdot }L^2(\mathbb{R};{\cal H})$ (by which we only mean
that $w\in C_+^0({\cal H})$ and that $\Vert w\Vert_{e^{\omega _1\cdot
  }L^2(\mathbb{R};{\cal H})}<\infty $, avoiding to define the larger
space $e^{\omega _1\cdot }L^2(\mathbb{R};{\cal H})$), then $Ew$
belongs to the same space and
\begin{eqnarray*}
\Vert Ew\Vert_{e^{\omega_1\cdot  }L^2(\mathbb{R};{\cal
    H})}&\le&\left(\int_0^\infty e^{-\omega _1t}\Vert S(t)\Vert
\,dt\right) \; \Vert w\Vert_{e^{\omega_1\cdot  }L^2(\mathbb{R};{\cal
    H})}\\
&\le & \frac{M}{\omega _1-\omega _0}\Vert w\Vert_{e^{\omega _1\cdot
  }L^2(\mathbb{R};{\cal H})}.
\end{eqnarray*}

\par Now we consider Laplace transforms. If $u\in e^{\omega \cdot
}{\cal S}(\mathbb{R};{\cal H})$, then the Laplace transform
$$
\widehat{u}(\tau )=\int_{-\infty}^{+\infty} e^{-t\tau }u(t)dt
$$
is well-defined in ${\cal S}(\Gamma _{\omega} ;{\cal H})$, where 
$$
\Gamma _{\omega} =\{ \tau \in \mathbb{C}; \Re \tau =\omega \}
$$
and we have Parseval's identity
\ekv{ny.3.5}
{
\frac{1}{2\pi }\Vert \widehat{u}\Vert_{L^2(\Gamma _{\omega} )}^2=\Vert
u\Vert_{e^{\omega \cdot }L^2}.
}

Now we make the assumptions in Theorem \ref{int2}, define $\omega $ and
$r(\omega )$ as there, and let $M,\omega _0$ be as above. Let $w\in e^{\omega
  \cdot }{\cal S}_+({\cal D}(A))$, where ${\cal S}_+({\cal D}(A))$ by
definition is the space of all $u\in {\cal S}(\mathbb{R}; {\cal
  D}(A))$, vanishing near $-\infty $. Then $w\in e^{\omega _1\cdot
}{\cal S}_+({\cal D}(A))$ for all $\omega _1\ge \omega $. If $\omega
_1>\omega _0$ then $u:=Ew$ belongs to $e^{\omega _1\cdot }{\cal
  S}_+({\cal D}(A))$ and solves (\ref{3}). Laplace transforming that
equation, we get 
\ekv{ny.4} { (\tau -A)\widehat{u}(\tau
  )=\widehat{w}(\tau ), } 
for $\Re \tau >\omega _0$. Notice here that
$\widehat{w}(\tau )$ is continuous in the half-plane $\Re \tau \ge
\omega $, holomorphic in $\Re \tau >\omega $, and
${{\widehat{w}}}_{\vert \Gamma _{\widetilde{\omega }}}\in {\cal
  S}(\Gamma _{\widetilde{\omega }})$ for every $\widetilde{\omega }\ge
\omega $. We use the assumption in the theorem to write
\ekv{ny.5}
{
\widehat{u}(\tau )=(\tau -A)^{-1}\widehat{w}(\tau ),
}
and to see that $\widehat{u}(\tau )$ can be extended to the half-plane
$\Re \tau \ge \omega $ with the same properties as $\widehat{w}(\tau
)$. By Laplace (Fourier) inversion from $\Gamma _{\omega } $ we conclude that $u\in
e^{\omega \cdot }{\cal S}_+({\cal D}(A))$. Moreover, since 
$$
\Vert \widehat{u}(\tau )\Vert_{{\cal H}}\le \frac{1}{r(\omega )}\Vert
\widehat{w}(\tau )\Vert_{{\cal H}},\ \tau \in \Gamma _{\omega} ,
$$
we get from Parseval's identity that 
\ekv{ny.6}
{
\Vert u\Vert_{e^{\omega \cdot }L^2}\le \frac{1}{r(\omega )} \Vert w\Vert_{e^{\omega \cdot }L^2}.
}

\par Using the density of ${\cal D}(A)$ in ${\cal H}$ together with
standard cutoff and regularization arguments, we see that
(\ref{ny.6}) extends to the case when $w\in e^{\omega \cdot
}L^2(\mathbb{R};{\cal H})\cap C_+^0({\cal H})$, leading to the fact
that $u:=Ew$ belongs to the same space and satisfies (\ref{ny.6}).

Consider $u(t)=S(t) v$, for $v\in D(A)$, solving the Cauchy problem
$$
\begin{array}{l}
(\pa_t -A) u=0\,,\, t\geq 0\,\,,\\
u(0)=v\,.
\end{array}
$$
Let $\chi$ be a decreasing Lipschitz function on $\mathbb R$, equal to
$1$ on $]-\infty,0]$ and vanishing near $+\infty$.
Then
$$
(\pa_t -A)(1- \chi ) u =-\chi'(t) u\,,
$$
and

\begin{eqnarray*}
\|\chi' u\|^2_{e^{\omega \cdot}L^2}
& =& \int_0^{+\infty}
 |\chi'(t)|^2 \|u(t)\| ^2  e^{-2 \omega t}\, dt\\
& 
\leq& \|\chi' m\|^2_{e^{\omega \cdot}L^2}\;\|v\|^2,
\end{eqnarray*}
where we notice that $\chi 'm$ is welldefined on $\mathbb{R}$ since
$\supp \chi '\subset [0,\infty [$.

Now $(1-\chi )u$, $\chi 'u$ are well-defined on $\mathbb{R}$, so 
\begin{equation}\label{13}
\|(1-\chi) u\|_{e^{\omega \cdot}L^2}
 \leq r(\omega)^{-1} \|\chi' u \|_{e^{\omega \cdot}L^2} \leq r(\omega)^{-1} 
\|\chi' m\|_{e^{\omega \cdot}L^2}\;\|v\|\,.
\end{equation}
Strictly speaking, in order to apply \eqref{ny.6}, we approximate $\chi$
 by a sequence of smooth functions.
Similarly, 
$$
\|\chi u\|_{e^{\omega \cdot}L^2(\mathbb R_+)}
 \leq \|\chi m\|_{e^{\omega \cdot}L^2(\mathbb{R}_+)}\;\|v\|\,,
$$
so
$$
\|u\|_{e^{\omega \cdot}L^2(\mathbb R_+)}
 \leq \left(r(\omega )^{-1}\|\chi' m\|_{e^{\omega \cdot}L^2} + \|\chi m\|_{e^{\omega
     \cdot}L^2(\mathbb{R}_+)}\right)
 \| v\|\,.
$$  

Let us now go from $L^2$ to $L^\infty$. For $t>0$, let 
$\chi_+(s)=\widetilde \chi (t-s)$ with $\widetilde \chi$ as $\chi$ above and
in addition
 $\supp \widetilde \chi \subset ]-\infty ,t]$, so that $\chi_+ (t) = 1$
 and $\supp\chi_+\subset [0,\infty [$.
Then
$$
\left(\pa_s -A\right) (\chi_+(s) u(s)) = \chi_+'(s) u(s)\,,
$$
and
$$
\chi_+ u (t) = \int_{-\infty}^{t} S(t-s)\, \chi_+'(s) \, u(s) \, ds\,.
$$

Hence, we obtain
\begin{equation}\label{16}
\begin{array}{ll}
e^{-\omega t} \| u(t)\|&= e^{-\omega t} \|\chi_+(t)  u(t)\|\\
& \leq \int_{-\infty}^t e^{-\omega t}\,m(t-s) |\widetilde \chi '(t-s)| \| u(s)\|\,ds\\
& \leq \int_{-\infty}^t e^{-\omega (t-s)}\,m(t-s)\, |\widetilde \chi '(t-s)|\;
e^{-\omega s} \| u(s)\|\,ds\\
&\leq \| m \widetilde \chi'\|_{e^{\omega \cdot}L^2}\; \| u\|_{e^{\omega
    \cdot}L^2(\supp\chi _+)}\,.
\end{array}
\end{equation}

Assume that
\begin{equation}\label{16a}
\chi =0 \hbox{ on } \supp \chi_+\, .
\end{equation}

Then $u$ can be replaced by $(1-\chi) u$ in the last line in
\eqref{16}
 and using \eqref{13} we get
\begin{equation}\label{17}
e^{-\omega t} \| u(t)\| \leq r(\omega )^{-1} \| m  \chi'\|_{e^{\omega \cdot}L^2}
\| m \widetilde \chi'\|_{e^{\omega \cdot}L^2} \|v\|\,.
\end{equation}
Let

\ekv{17a}{
\supp \chi \subset ]-\infty ,a]\,,\, \supp \widetilde \chi \subset
]-\infty ,\widetilde a ],
a+\widetilde a = t\, ,}
so that (\ref{16a}) holds.

For a given $a>0$, we look for $\chi$ in (\ref{17a}) such that $\| m
\chi'\|_{e^{\omega \cdot}L^2}$
 is as small as possible. By the Cauchy-Schwarz inequality, 
\begin{equation}\label{18}
1 = \int_0^a |\chi'(s) | ds \leq \|\chi' m\|_{e^{\omega \cdot}L^2}
\|\frac 1m\|_{e^{-\omega \cdot}L^2(]0,a[)}\,,
\end{equation}
so
\begin{equation}\label{19}
 \|\chi' m\|_{e^{\omega \cdot}L^2}  \geq \frac{1}{\|\frac
   1m\|_{e^{-\omega \cdot}L^2(]0,a[)}}\,.
\end{equation}
We get equality in \eqref{19} if for some constant $C$,
$$
|\chi'(s)| m(s) e^{-\omega s} = C \frac{1}{m(s)} e^{\omega s}, \hbox{
  on }[0,a],
$$
i.e.
$$
\chi'(s) m(s) e^{-\omega s} = - C \frac{1}{m(s)} e^{\omega s}, \hbox{
  on }[0,a],
$$
where $C$ is determined  by the condition $1 = \int_0^a |\chi'(s) | ds$.\\
We get
$$
C= \frac{1}{\|\frac
   1m\|^2_{e^{-\omega \cdot}L^2(]0,a[)}}\,,
$$
Here $\chi (s)=1$ for $s\le 0$,  $\chi (s)=0$ for $s\ge a$,
 $$
\chi (s)=C\int_s^a
\frac{1}{m(\sigma )^2}e^{2\omega \sigma }d\sigma ,\ 0\le s\le a. $$

\par With the similar optimal choice of $\widetilde{\chi }$, for which
$$
\Vert \widetilde{\chi }'m\Vert_{e^{\omega \cdot }L^2}=\frac{1}{\Vert
  \frac{1}{m}\Vert_{e^{-\omega \cdot }L^2([0,\widetilde{a}])}},
$$ 
we get from (\ref{17}): \ekv{20} { e^{-\omega t}\Vert u(t)\Vert \le
  \frac{\Vert v\Vert}{r(\omega )\Vert \frac{1}{m}\Vert_{e^{-\omega \cdot
      }L^2([0,a])}\Vert \frac{1}{m}\Vert_{e^{-\omega \cdot
      }L^2([0,\widetilde{a}])}}, } provided that
$a,\widetilde{a}>0$, $a+\widetilde{a}=t$, for any $v\in
D(A)$. Observing that $D(A)$ is dense in $\mathcal H$,  this completes the proof of
Theorem \ref{int2}.
\subsection{Proof of Theorem \ref{int3}}\label{pr3} 
We can apply Theorem \ref{int2} to the restriction $\widetilde{S}(t)$
of $S(t)$ to the range \break  ${\cal R}(1-\Pi _+)$ of $1-\Pi_+$. The
generator is the restriction $\widetilde{A}$ of $A$ so we get 
\ekv{int3.1} { \Vert \widetilde{S}(t)\Vert \le
  \frac{e^{\widetilde{\omega} t}}
{r(\widetilde{\omega} )\Vert \frac{1}{m}\Vert_{e^{-\widetilde{\omega} \cdot
      }L^2([0,a])}\Vert \frac{1}{m}\Vert_{e^{-\widetilde{\omega} \cdot
      }L^2([0,\widetilde{a}])}}.}
Then (\ref{int.5}) follows from the fact that
$R(t)=\widetilde{S}(t)(1-\Pi _+)$.
\appendix 

\section{An iterative improvement of  Theorem \ref{int2}}

Working entirely on the semi-group side and applying Theorem
\ref{int2} repeatedly, we shall see how to gain an
extra decay ${\cal O}(1)\exp (-t^{1/2}/C)$ for some $C>0$. It is not
clear that this result is of practical use, especially in view of
Lemma \ref{int01}, but the computations are amusing.

\medskip Recall that under the assumptions of Theorem \ref{int2} we have the
estimate \eqref{int.4}. Here we may have $m$ bounded continuous for $0\le t<T$
and equal to $+\infty $ for $t\ge T$, where $T>0$.

Write $m(t)=\widetilde{m}(t)e^{\omega t}$. Then \eqref{int.4} shows that $\Vert
S(t)\Vert \le \widehat{m}(t)e^{\omega t}$, where 
\begin{equation}\label{rem.1}
\widehat{m}(t)\le \frac{1}{r(\omega )\Vert
  \frac{1}{\widetilde{m}}\Vert_{[0,a]}\Vert
  \frac{1}{\widetilde{m}}\Vert_{[0,\widetilde{a}]}}, \ a+\widetilde{a}=t.
\end{equation}
Take $a=\widetilde{a}=t/2$ and divide the previous inequality by
$r(\omega )$:
$$
\frac{\widehat{m}(t)}{r(\omega ) }\le
 \frac{1}{\int_0^{t/2}(\frac{r(\omega)}{\widetilde{m}(s)})^2\, ds},
$$
which we can also write 
$$
\widehat{f}(t)\ge \int_0^{t/2}\widetilde{f}(s)^2 \,ds\,,\
\widetilde{f}(t):=\frac{r(\omega)}{\widetilde{m}(t)},\, 
\widehat{f}(t):=\frac{r(\omega)}{\widehat{m}(t)}\,.
$$

\par Now assume that $e^{-\omega t}\Vert S(t)\Vert \le
\widetilde{m}(t)\le {\cal O}(1)$ for $0\le t<T$. Then we extend
$\widetilde{m}$ to $[0,+\infty [$, by defining
\begin{equation}\label{rem.2}
\frac{\widetilde{m}(t)}{r(\omega)}=\frac{1}{\int_0^{t/2}(\frac{r(\omega)}{\widetilde{m}(s)})^2\,
ds},
\end{equation}
first for $T\le t<2T$, then for $2T\le t<4T$ and so on. 
Correspondingly, we have
\begin{equation}\label{rem.3}
\widetilde{f}(t)= \int_0^{t/2}\widetilde{f}(s)^2 \, ds,\ t\ge T.
\end{equation}

Theorem \ref{int2}
now shows that $e^{-\omega t}\Vert S(t)\Vert \le
\widetilde{m}(t)\le {\cal O}(1)$ for all $t\ge 0$. By construction we
see that $\widetilde{m}(t)$ is decreasing on $[T,+\infty [$, so we
have 
\begin{equation}\label{rem.4}
e^{-\omega t}\Vert S(t)\Vert \le M,\ M=\max(\sup_{[0,T[}\widetilde{m},
\frac{1}{r(\omega )\int_0^{T/2}\widetilde{m}(s)^{-2}\,ds}). 
\end{equation}

\par Notice that $\widetilde{f}$ is increasing on $[T,+\infty [$. We
look for upper bounds on $\widetilde{m}$ or equivalently for lower
bounds on $\widetilde{f}$. For $k\ge 1$, put $I_k=[T2^{k-1},T2^k[$, so
that the length of $I_k$ is $|I_k|=T2^{k-1}$. Put 
$$
F(k)=\inf_{I_k}\widetilde{f}=\widetilde{f}(T2^{k-1}) \hbox{ when }k\ge
1,\quad F(0)=\inf_{[0,T[}\widetilde{f}(t)\,.
$$
Then,
$F(1)=\int_0^{T/2}\widetilde{f}(t)^2\,dt\ge \frac{T}{2}F(0)^2$, which we write
$$TF(1)\ge \frac{1}{2}(TF(0))^2\, .$$
For $k\ge 1$, we get 
$$
F(k+1)\ge \int_0^{T2^{k-1}}\widetilde{f}(t)^2dt\ge TF(0)^2+TF(1)^2+2TF(2)^2+...+2^{k-2}TF(k-1)^2,
$$
which we write
\begin{equation}\label{rem.5}
TF(k+1)\ge (TF(0))^2+(TF(1))^2+2(TF(2))^2+...+2^{k-2}(TF(k-1))^2.
\end{equation}
Since $\widetilde{f}$ is increasing on $[T,+\infty [$, we have 
$$
F(1)\le F(2)\le F(3)\le ...
$$
Thus for $k\ge 2$,
$$TF(k+1)\ge 2^{k-2}(TF(k-1))^2\ge 2^{k-2}(TF(1))^2\ge
2^{k-4}(TF(0))^4\,,$$ which we write
$$
TF(k)\ge 2^{k-5}(TF(0))^4\,,\ k\ge 3.
$$
Let $k_0$ be the smallest integer $k\ge 3$ such that 
$$
2^{k-5}(TF(0))^4\ge 2,
$$
so that $TF(k)\ge 2$ for $k\ge k_0\,$.

\par Now return to (\ref{rem.5}) which implies that 
$$TF(k+1)\ge 2^{k-2}(TF(k-1))^2,\ k\ge 1.$$
We get 
$$TF(k+2)\ge 2^{k-1}(TF(k))^2,\ k\ge 1,$$
implying,
$$
T(F(k+2))\ge (TF(k))^2,\quad \ln (TF(k+2))\ge 2\ln (TF(k)).
$$
In particular,
$$\ln (TF(k_0+2\nu ))\ge 2^{\nu }\ln (TF(k_0))\ge 2^{\nu }\ln 2,\ \nu
\in \mathbb{N}.$$
We conclude that 
$$
T\widetilde{f}(t)\ge 2^{2^\nu },\ 2^{k_0+2\nu -1}\le t/T<2^{k_0+2\nu }.
$$
The last inequality for $t$
 implies that $2^\nu >(2^{-k_0}t/T)^{1/2}$, 
so we get 
\begin{equation}\label{rem.6}
T\widetilde{f}(t)\ge 2^{(2^{-k_0}t/T)^{1/2}},\ t/T\ge 2^{k_0-1},
\end{equation}
or equivalently,
\begin{equation}\label{rem.7}
\frac{\widetilde{m}(t)}{r(\omega )T}\le 2^{-(2^{-k_0}t/T)^{1/2}},\ t/T\ge 2^{k_0-1},
\end{equation}
where we recall that $k_0$ is the smallest integer such that 
\begin{equation}\label{rem.8}
2^{k_0}\ge \max (\frac{2^6}{(TF(0))^4},8).
\end{equation}


\newpage 
\begin{thebibliography}{99}

\bibitem{Alm} Y. Almog.
\newblock The stability of the normal state of superconductors
 in the presence of electric currents.
\newblock Siam J. Math. Anal. 40 (2) (2008), p.~824-850.

\bibitem{Bou02}
L.S.~Boulton. \newblock Non-self-adjoint harmonic oscillator, compact
semigroups and pseudospectra. \newblock J. Operator Theory  47(2) (2002), p.~413-429.

\bibitem{BuZw} N.~Burq, M.~Zworski, \newblock
  Geometric control in the presence of a black box. \newblock
  J. Amer. Math. Soc.  17(2) (2004), p.~443-471. 

\bibitem{Dav2} E.B.~Davies.
 \newblock Semi-classical states for non self-adjoint Schr\"odinger operators.
 \newblock Comm. in Math. Phsics 200  (1999), p.~35-41.

\bibitem{Dav3} E.B.~Davies, \newblock Pseudospectra, the harmonic
  oscillator and complex resonances. \newblock Proc. Roy. Soc. 
London Ser. A 455 (1999), p.~585-599.

\bibitem{Da07} E.B.~Davies.  \newblock {\it Linear operators and their spectra,}
   \newblock Cambridge Studies in Advanced Mathematics, 106. Cambridge University
  Press, Cambridge, 2007.


\bibitem {DSZ} N.~Dencker, J.~Sj\"ostrand, and M.~Zworski.
\newblock Pseudospectra of semi-classical (pseudo)differential
operators.
\newblock Comm. Pure Appl. Math. 57 (3) (2004), p.~384-415.

\bibitem{EnNa07} K.J.~Engel, R.~Nagel.  \newblock {\it One-parameter semigroups
  for linear evolution equations,}  \newblock Graduate
  Texts in Mathematics, 194. Springer-Verlag, New York, 2000.

\bibitem{EnNa2} K.J.~Engel, R.~Nagel.
 \newblock {\it A short course on operator semi-groups},  \newblock  Unitext, Springer-Verlag (2005).


\bibitem{GGN} I. Gallagher, T. Gallay and F. Nier.
\newblock  Spectral asymptotics for large skew-symmetric perturbations of
the harmonic oscillator,
\newblock  Preprint 2008.

\bibitem{He} B. Helffer. 
 \newblock On spectral problems related to a time dependent 
 model in superconductivity
 with electric current.  \newblock Proceedings of the Colloque
  sur les \'equations aux d\'eriv\'ees partielles, \'Evian, June 2009,
  to appear.


\bibitem {HelNi} B.~Helffer and F.~Nier.
 \newblock {\it Hypoelliptic estimates and  spectral theory
 for  Fokker-Planck operators 
and Witten Laplacians,}
 \newblock no.~1862 in  Lecture Notes in Mathematics, Springer-Verlag (2004).

\bibitem{HeHiSj08a} F.~H{\'e}rau, M.~Hitrik, J.~Sj\"ostrand.
  \newblock  Tunnel effect for Fokker-Planck type operators. \newblock
Annales Henri Poincar\'e, 9(2) (2008), p.~209-274.

\bibitem{HeHiSj08b}F.~H{\'e}rau, M.~Hitrik, J.~Sj\"ostrand.  \newblock
  Tunnel effect for
Kramers-Fokker-Planck type operators: return to equilibrium and
applications.  \newblock 
International Math Res Notices, Vol. 2008, Article ID rnn057, 48p.

\bibitem{HeNi04} F.~H\'erau, F.~Nier.
  \newblock Isotropic hypoellipticity and trend to equilibrium for the
  Fokker-Planck equation with a high-degree potential.A \newblock rch. Ration. Mech. Anal. 171(2) (2004), p.~151-218.


\bibitem{HeSjSt05}F.~H{\'e}rau, J.~Sj\"ostrand, C.~Stolk.
  \newblock Semiclassical analysis for the Kramers-Fokker-Planck equation. 
   \newblock Comm. PDE 30(5--6) (2005), p.~689-760.

\bibitem{Hi} M. Hitrik.
 \newblock Eigenfunctions and expansions for damped wave equations.
 \newblock Meth. Appl. Anal. 10 (4) (2003), p.~1-22.

\bibitem{Paz} A. Pazy.
 \newblock {\it Semigroups of linear operators and applications to partial
  differential operators.}
 \newblock Appl. Math. Sci. Vol. 44, Springer (1983).

\bibitem{Pr06} K.~Pravda-Starov. \newblock A complete study of the
  pseudo-spectrum for the rotated harmonic oscillator. \newblock
  J. London Math. Soc. (2)  73(3) (2006), p.~745-761.

\bibitem{Sch09} E.~Schenk,  {\it Syst\`emes quantiques ouverts et
    m\'ethodes semi-classiques,} th\`ese novembre 2009.\\
http://www.lpthe.jussieu.fr/~schenck/thesis.pdf

\bibitem {Sj} J. Sj\"ostrand.  \newblock Resolvent estimates  for non-self-adjoint operators via
semi-groups.  \newblock To appear in a volume in honor of V.~Maz'ya,\\
http://arxiv.org/abs/0906.0094.

\bibitem{Sj2} J. Sj\"ostrand.  \newblock  Spectral properties for non
    self-adjoint differential operators.   \newblock Proceedings of the Colloque
  sur les \'equations aux d\'eriv\'ees partielles, \'Evian, June 2009,
  to appear.


\bibitem {Tr}  L.N. Trefethen. \newblock  Pseudospectra of linear
  operators.  \newblock 
SIAM Rev. 39 (3) (1997), p.~383-406.

\bibitem {TrEm} L.N.~Trefethen, M.~Embree.  \newblock {\it Spectra and
  pseudospectra. The behavior of nonnormal matrices and operators.}
  \newblock  Princeton University Press, Princeton, NJ, 2005.


\bibitem{Vi} C. Villani.  \newblock {\em Hypocoercivity.}  \newblock Memoirs of the AMS,
  Vol. 202, no. 950 (2009).

\bibitem{Zw} M.~Zworski.
 \newblock A remark on a paper by E.B.~Davies.  \newblock 
 Proc. Amer. Math. Soc. 129  (2001),
p.~2955-2957.


\end{thebibliography}
\end{document}